\documentclass[12pt,twoside]{article}
\usepackage{amsmath, amsthm, amssymb}
\usepackage{color}
\usepackage{hyperref}

\theoremstyle{plain}
\swapnumbers

\begin{document}
\font\psaci=rsfs10 

\newcommand{\R}{\mathbb{R}}
\newcommand{\N}{\mathbb{N}}
\newcommand{\dd}{\hskip0.2mm\mbox{\rm d}}
\newcommand{\wt}{\widetilde}
\newcommand{\wh}{\widehat}

\newcommand{\inn}{\,{\in}\,}
\newcommand{\ab}{[a,b\,]}
\newcommand{\intab}{\int_a^b}

\newcommand{\var}{\mbox{\rm var}}
\newcommand{\varab}{\var_a^b}

\newcommand{\dis}{\displaystyle}

\newcommand{\BV}{{\rm{BV}}}
\newcommand{\G}{{\rm{G}}}

\newcommand{\eps}{\varepsilon}
\newcommand{\ph}{\varphi}
\newcommand{\balfa}{{\boldsymbol{\alpha}}}
\newcommand{\bbeta}{{\boldsymbol{\beta}}}
\newcommand{\ksi}{{\boldsymbol{\xi}}}
\newcommand{\brho}{{\boldsymbol{\rho}}}
\newcommand{\bsigma}{{\boldsymbol{\sigma}}}
\newcommand{\eeta}{{\boldsymbol{\eta}}}
\newcommand{\dzeta}{{\boldsymbol{\zeta}}}

\newcommand{\opD}{{\mbox{\psaci D}\hskip0.75mm}}
\newcommand{\opG}{{\mbox{\psaci G}\hskip0.75mm}}

\def\skipaline{\removelastskip\vskip10pt plus 1pt minus 1pt}
\def\skiphalfaline{\removelastskip\vskip5pt plus 1pt minus 1pt}

\newcommand{\et}{\mbox{ \ and \ }}
\newcommand{\ia}{\noalign{\noindent\mbox{and}}}

\newtheorem{The}{\,Theorem}[section]
\newtheorem{Lem}[The]{\,Lemma}
\newtheorem{Pro}[The]{\,Proposition}
\newtheorem{Cor}[The]{\,Corollary}
\newtheorem{Cors}[The]{\,Corollaries}

\theoremstyle{definition}
\newtheorem{Not}[The]{\,Notation}
\newtheorem{Rem}[The]{\,Remark}
\newtheorem{Rems}[The]{\,Remarks}
\newtheorem{Cvi}[The]{\,Exercise}
\newtheorem{Cvs}[The]{\,Exercises}
\newtheorem{Exa}[The]{\,Example}
\newtheorem{Exs}[The]{\,Examples}
\newtheorem{Def}[The]{\,Definition}

\numberwithin{equation}{section}
\numberwithin{The}{section}

\parindent 4mm
\pagestyle{empty}

\title{\large\bf On the relationships between Stieltjes type integrals
\\ of Young, Dushnik and Kurzweil}
\author{
{\normalsize Umi Mahnuna Hanung, Yogyakarta}\\[-1mm]
{\small e-mail: {\tt{hanungum@\allowbreak ugm.ac.id}}}
\\
{\normalsize Milan Tvrd\'y, Praha}\\[-1mm]
{\small e-mail: {\tt{tvrdy@\allowbreak math.cas.cz}}}
}
\date{\today}

\maketitle

Integral equations of the form
\[
\quad x(t)=x(t_0)+\int_{t_0}^t\dd[A]\,x=f(t)-f(t_0)
\]
are natural generalizations of systems of linear differential equations. Their
main goal is that they admit solutions which need not be absolutely continuous.
Up to now such equations have been considered by several authors starting with
J.~Kurzweil \cite{Ku} and T.H.~Hildebrandt \cite{Hi0}. For further contributions
see e.g. \cite{As}, \cite{hoe}, \cite{MeZh}, \cite{MST}, \cite{ssy}--\cite{stv}
and references therein. These papers worked with several different concepts of
the Stieltjes type integral like Young's (Hildebrandt), Kurzweil's (Kurzweil,
Schwabik and Tvrd\'{y}),  Dushnik's (H\"{o}nig) or Lebesgue's (Ashordia, Meng
and Zhang). Thus an interesting question arises: {\bf what are the relationships
between all these concepts}?

It is known that (cf. \cite[Theorem 1.2.1]{Ku}) the Kurzweil-Stieltjes integral
is in finite dimensional setting equivalent with the Ward-Perron-Stieltjes, while
the relationship between the Ward-Perron-Stieltjes and the Lebesgue-Stieltjes integrals
has been described in \cite[Theorem VI.8.1]{sa}. For more details, see Chapter 6
of \cite{MST}. The relationship between the Young and the Dushnik integrals is
indicated by \cite[Theorem B]{mcn}. Finally, the relationship between the Young
integral and the Kurzweil-Stieltjes one has been considered in \cite{ssy} and
\cite{ssm}. Our aim is to complete this schedule. In addition, we will present
also convergence results that are possibly new for the Young and Dushnik integrals.
Let us emphasize that the proofs of all the assertions  presented in this paper are
based on rather elementary tools.

\section{Preliminaries}
In this paper the symbols like $\R,$ $\N,$ $\ab,$ $(a,b),$ $\varab\,f$ and $\|f\|_\infty$
have their usual and traditional meaning. Furthermore, recall that a finite sequence
$\balfa=\{\alpha_0,\dots,\alpha_{\nu(\balfa)}\}$ of points from $\ab$ is a division of $\ab$ if
$a=\alpha_0<\dots<\alpha_{\nu(\balfa)}=b.$ The couple $P\,{=}\,(\balfa,\ksi)$ is a partition
of $\ab$ if $\balfa$ is a division of $\ab$ and $\ksi\,{=}\,\{\xi_1,\dots,\xi_{\nu(P)}\}$
is a finite sequence such that $\xi_j\in[\alpha_{j-1},\alpha_j]$ for all $j.$ If
$P\,{=}\,(\balfa,\ksi)$ is a partition of $\ab,$ the elements of $\balfa$ and $\ksi$ are
always denoted respectively as $\alpha_j$ and $\xi_j.$ At the same time the number of elements
of $\ksi$ is denoted by $\nu(P).$ (If $P=(\balfa,\ksi),$ then $\nu(P)=\nu(\balfa).$)

Recall that a function $f:\ab\to\R$ is regulated on $\ab$ if it has finite one sided limits
\[
   \lim_{\tau\to t-}f(\tau)=f(t-) \et \lim_{\tau\to s+}f(\tau)=f(s+)
\]
for all $t\in(a,b]$ and $s\in[a,b).$ For every function $f$ regulated on $\ab$ and points
$t\in(a,b]$ and $s\in[a,b),$ we denote
\[
   \Delta^-f(t)=f(t)-f(t-) \et \Delta^+f(s)=f(s+)-f(s).
\]
The set of all functions regulated on $\ab$ is denoted by $\G(\ab).$

Furthermore, a function $f:\ab\to\R$ is a finite step function if there exist an $m\in\N,$
sequences $\big\{\wt{c}_k:k\in\{1,...,m\}\big\}\subset\R,$
$\big\{\wt{d}_k:k\in\{1,...,m\}\big\}\subset\R,$ and a division
$\bsigma=\{\sigma_0,\dots,\sigma_m\}$ of $\ab$ such that $f(\sigma_k)=\wt{c}_k$ for
$k\in\{0,\dots,m\}$ and $f(x)=\wt{d}_k$ for $k\in\{1,\dots,m\}$ and
$x\in(\sigma_{k-1},\sigma_k),$ i.e.
\[
   f(x)=\sum_{k=0}^m \wt{c}_k\,\chi_{[\sigma_k]}(x)
        +\sum_{k=1}^m \wt{d}_k\,\chi_{(\sigma_{k{-}1},\sigma_k)}(x).
\]
Equivalently,
\begin{equation}\label{step.funct}
   f(x)=c+\sum_{k=0}^{m-1}c_k\,\chi_{(\sigma_k,b]}(x)
         +\sum_{k=1}^{m-1}d_k\,\chi_{[\sigma_k,b]}(x)+d\,\chi_{[b]}(x)
   \quad\mbox{for \ } x\in\ab,
\end{equation}
where $c=\wt{c}_0,\, c_k=\wt{d}_{k{+}1}-\wt{c}_k$ for $k\in\{0,\dots,m\},$
$d_k=\wt{c}_k-\wt{d}_k$ for $k\in\{1,\dots,m-1\}$ and $d=\wt{c}_m-\wt{d}_m.$
Then,
\begin{gather*}
   f(a)=c,\,\,f(x-)=f(x) \mbox{\ for\ } x\in(a,b]\setminus\{\sigma_k\},\,\,
              f(x+)=f(x) \mbox{\ for\ } x\in[a,b)\setminus\{\sigma_k\},
\\\ia
   \Delta^+f(\sigma_k)=c_k \mbox{ \ for\ } k\in\{0,\dots,m-1\},\quad
   \Delta^-f(\sigma_k)=d_k \mbox{ \ for\ } k\in\{1,\dots,m\}.
\end{gather*}

\skiphalfaline

It is known (cf. e.g. \cite[Theorem 3.1]{hoe1}) that $f:\ab\to\R$ is regulated
if and only if it is a uniform limit of step functions.

Further recall that for a given function $f:\ab\to\R$ the symbol $\varab\,f$
stands for its variation and $\|f\|=\dis\sup_{t\in\ab}|f(t)|.$

\skipaline

For functions $f,g{:}\,\ab{\to}\R$ and a partition $P\,{=}\,(\balfa,\ksi)$ of
$\ab$ we set
\begin{align*}
   S(f,\dd g,P)&{=}\sum_{j=1}^{\nu(P)}f(\xi_j)\,[g(\alpha_j)-g(\alpha_{j-1})]
  \\\noalign{\noindent\mbox{and, if $g$ is regulated,}}
   S_Y(f,\dd g,P)&{=}\sum_{j=1}^{\nu(P)}
            \Big(f(\alpha_{j-1})\,\Delta^+g(\alpha_{j{-}1})
               \,{+}\,f(\xi_j)\,[g(\alpha_j-){-}g(\alpha_{j{-}1}+)]
               \,{+}\,f(\alpha_j)\,\Delta^-g(\alpha_j)\Big)
\end{align*}
and define:
\begin{itemize}
\item{} The Young integral $({\rm Y})\int_a^b f\,\dd g$ \
  \big(the Dushnik integral $({\rm D})\int_a^b f\,\dd g$\big)
  exists and equals $I\,{\in}\,\R$ if
\vskip-3mm
\[
\begin{array}{l}
   \mbox{for every\ } \eps>0\,\,\mbox{there is a division \ } \balfa_\eps
   \mbox{ \ of $\ab$ \ such that}
  \\[1mm]\hskip10mm
   |S_Y(f,\dd g,P)-I|<\eps \quad\big(\mbox{or}\quad |S(f,\dd g,P)-I|<\eps\big)
  \\[1mm]
   \mbox{\ holds for all partitions \ }
   P=(\balfa,\ksi)\mbox{\ of $\ab$ such that\ } \balfa\supset\balfa_\eps
   \mbox{ \ and}
  \\[1mm]\hskip10mm
   \alpha_{j-1}<\xi_j<\alpha_j \quad\mbox{for all \ } j\in\{1,\dots,\nu(\balfa)\}.
\end{array}
\]

\item{} The Kurzweil-Stieltjes integral $({\rm K})\int_a^b f\,\dd g$ exists and equals
$I\,{\in}\,\R$ if
\vskip-7mm
\[
\begin{array}{l}
  \mbox{for every \ } \eps>0 \mbox{ \ there exists a function \ $\delta_\eps:\ab\to(0,1)$
  \ such that}
 \\[1mm]\hskip12mm
  |I-S(f,\dd g,P)|<\eps
 \\[1mm]
  \mbox{holds for all partitions $P=(\balfa,\ksi)$ of $\ab$ such that}
 \\[1mm]\hskip12mm
  [\alpha_{j-1},\alpha_j]\subset[\xi_j-\delta_\eps(\xi_j),\xi_j+\delta_\eps(\xi_j)].
\end{array}
\]
\end{itemize}

\section{Results}

Our main goal is the following assertion:

\skiphalfaline

\begin{The}\label{main}
Suppose that $f$ and $g$ are regulated on~$\ab$ and at least one of them has a bounded
variation on $\ab.$ Then all the integrals ${\rm(K)}\intab f\,\dd g,$ ${\rm(Y)}\intab f\,\dd g$
and ${\rm(D)}\intab f\,\dd g$ exist and
\begin{equation}\label{eq-main}
    ({\rm K})\!\intab f\,\dd g=({\rm Y})\!\intab f\,\dd g
    =f(b)\,g(b)-f(a)\,g(a)-({\rm D})\!\intab g\,\dd f.
\end{equation}
\end{The}

\skipaline

To prove Theorem \ref{main}, we will need several auxiliary results. First, we
we will restrict ourselves to some simpler special cases.

\begin{Lem}\label{L1}
\begin{itemize}
\item[{\rm \,\,(i)}]
The equalities \eqref{eq-main} hold for every $f\,{:}\,\ab\,{\to}\,\R$ whenever $g$
is a finite step function.
\item[{\rm (ii)}]
The equalities \eqref{eq-main} hold for every $g\,{:}\,\ab\,{\to}\,\R$ is regulated
and $f$ is a finite step function.
\end{itemize}
\end{Lem}
\begin{proof}
Let $g\in\G(\ab),$ $\tau\in(a,b)$ and let the functions $f_i,$ $i\in\{1,\dots,5\},$ be
defined on $\ab$ by
\[
  f_1=1,\quad  f_2=\chi_{(\tau,b]},\quad f_3=\chi_{(\tau,b]},\quad f_4=\chi_{[\tau,b]},
  \quad f_5:={\chi}_{[b]}.
\]
Then
\begin{equation}\label{r1}
\left.\begin{array}{l}\dis
    ({\rm K})\!\intab f_1\,\dd g=({\rm Y})\!\intab f_1\,\dd g
    =({\rm D})\!\intab f_1\,\dd g=g(b)-g(a),
\\[4mm]\dis
   ({\rm K})\!\intab f_2\,\dd g=({\rm Y})\!\intab f_2\,\dd g=g(b)-g(a+),\,
   ({\rm D})\!\intab f_2\,\dd g=g(b)-g(a),
\\[4mm]\dis
   ({\rm K})\!\intab f_3\,\dd g=({\rm Y})\!\intab f_3\,\dd g=g(b)-g(\tau+),\,
   ({\rm D})\!\intab f_3\,\dd g=g(b)-g(\tau),
\\[4mm]\dis
   ({\rm K})\!\intab f_4\,\dd g=({\rm Y})\!\intab f_4\,\dd g=g(b)-g(\tau-),\,
   ({\rm D})\!\intab f_4\,\dd g=g(b)-g(\tau)
\\[4mm]\dis
   ({\rm K})\!\intab f_5\,\dd g=({\rm Y})\!\intab f_5\,\dd g=\Delta^-f(b),\,
   ({\rm D})\!\intab f_5\,\dd g=0
\end{array}\right\}
\end{equation}
and
\begin{equation}\label{r2}
\left.\begin{array}{l}\dis
    ({\rm K})\!\intab g\,\dd f_1=({\rm Y})\!\intab g\,\dd f_1
    =({\rm D})\!\intab  g\,\dd f_1=0,
\\[4mm]\dis
   ({\rm K})\!\intab g\,\dd f_2=({\rm Y})\!\intab g\,\dd f_2=g(a),\,
   ({\rm D})\!\intab g\,\dd f_2=g(a+),
\\[4mm]\dis
   ({\rm K})\!\intab g\,\dd f_3=({\rm Y})\!\intab g\,\dd f_3=g(\tau),\,
   ({\rm D})\!\intab g\,\dd f_3=g(\tau+),
\\[4mm]\dis
   ({\rm K})\!\intab g\,\dd f_4=({\rm Y})\!\intab g\,\dd f_4=f(\tau),\,
   ({\rm D})\!\intab g\,\dd f_4=g(\tau-),
\\[4mm]\dis
   ({\rm K})\!\intab g\,\dd f_5=({\rm Y})\!\intab g\,\dd f_5=f(b),\,
   ({\rm D})\!\intab g\,\dd f_5=g(b-).
\end{array}\hskip25mm\right\}
\end{equation}
Indeed, given an arbitrary partition $P=(\balfa,\ksi)$ of $\ab$ such that
$\xi_j\in(\alpha_{j-1},\alpha_j)$ for all $j\,{\in}\,\{1,\dots,\nu(P)\},$
we get
\[
   S_Y(f_3,\dd g,P)=g(b)-g(\tau+) \et S_D(f_3,\dd g,P)=g(b)-g(\tau),
\]
wherefrom the equalities $({\rm Y})\!\intab f_2\,\dd g=g(b)-g(\tau+)$ and
$({\rm D})\!\intab f_3\,\dd g=g(b)-g(\tau)$ from \eqref{r1} immediately
follow. Similarly, we can justify all the other relations given
in \eqref{r1}--\eqref{r2}.

\skiphalfaline

Now, having in mind \eqref{r1}--\eqref{r2} we can deduce that
\[
    ({\rm K})\!\intab f_i\,\dd g=({\rm Y})\!\intab f_i\,\dd g
    =f_i(b)\,g_i(b)-f_i(a)\,g_i(a)-({\rm D})\!\intab g\,\dd f_i.
\]
holds for every $i\in\{1,\dots,5\}.$ Finally, since by \eqref{step.funct} every finite
step function is a linear combination of functions of the type $\{f_1,\dots,f_5\},$
the proof of the lemma easily follows.
\end{proof}

\skipaline

Estimates needed later are summarized in the following lemma.

\begin{Lem}\label{L2}
Let $g\in\G(\ab),$ $f:\ab\to\R$ and a partition $P$ of $\ab$ be given. Then
the estimates
\begin{align}\label{sum-1}
&\left.\begin{array}{l}
  \big|\,\,S(f,\dd g,P)\big|\le\|f\|_\infty\,\varab g,
 \\[4mm]
  \big|\,\,S(f,\dd g,P)\big|\le(|f(a)|+|f(b)|+\varab f)\,\|g\|_\infty,
\end{array}\hskip32mm\right\}
 \\\ia\label{sum-2}
&\left.\begin{array}{l}
  \big|S_Y(f,\dd g,P)\big|\le\|f\|_\infty\,\varab g,
 \\[4mm]
  \big|S_Y(f,\dd g,P)\big|\le(|f(a)|+|f(b)|+\varab f)\,\|g\|_\infty
\end{array}\hskip32mm\right\}
\end{align}
are true. Furthermore, the estimates
\begin{align}\label{est-1}
  &\hskip11mm\Big|\intab f\,\dd g\Big|\le\|f\|_\infty\,\varab g
\\\ia\label{est-2}
  &\hskip11mm\Big|\intab f\,\dd g\Big|\le(|f(a)|+|f(b)|+\varab f)\,\|g\|_\infty
\end{align}
hold for each of the three integrals under consideration, whenever it exist.
\end{Lem}
\begin{proof}
For the Kurzweil-Stieltjes integral these inequalities are well-known, cf. e.g. \cite{tv}
or Chapter 6 of \cite{MST}. Since the set of admissible partitions for the Dushnik integral
is contained in that for the Kurzweil-Stieltjes integral, it follows immediately that
relations \eqref{sum-1}, \eqref{est-1} and \eqref{est-2} are true also for the Dushnik integral.
So, it remains to consider the Young integral.

\skiphalfaline

\noindent
a) \ If $a\le\alpha\le\xi\le\beta\le b,$ then
\begin{align*}
  &|f(\alpha)\,\Delta^-g(\alpha)+f(\xi)\,[g(\beta-)-g(\alpha+)]+f(\beta)\,\Delta^-g(\beta)|
\\[1mm]
  &\quad\le\|f\|_\infty\,\big(|\Delta^-g(\alpha)|+|g(\beta-)-g(\alpha+)|+|\Delta^-g(\beta)|)
        \le\|f\|_\infty\,\var_\alpha^\beta\,g,
\end{align*}
wherefrom it is easy to deduce that the estimate
\[
  |S_Y(f,\dd g,P)\le\|f\|_\infty\,\varab g
\]
holds for every partition $P$ of $\ab$  This means that the former inequality from \eqref{sum-2}
and estimate \eqref{est-1} are true also for the Young integral.

\skiphalfaline

\noindent
b) \ Observe that
\begin{align*}
   &\hskip-8mm
    f(\alpha)\,[g(\alpha+)-g(\alpha)]+f(\xi)\,[g(\beta-)-g(\alpha+)]
   +f(\beta)\,[g(\beta)-g(\beta-)]
  \\[1mm]&\hskip-5mm
  =[f(\alpha){-}f(\xi)]\,g(\alpha+)+[f(\xi)-f(\beta)]\,g(\beta-)
    +f(\beta)\,g(\beta)-f(\alpha)\,g(\alpha)
\end{align*}
holds for all $\alpha,\xi,\beta\in\ab$ such that $a\le\alpha\le\xi\le\beta\le b.$
Having this in mind we can verify the estimate
\[
  |S_Y(f,\dd g,P)\le(|f(a)|+|f(b)|+\varab f)\,\|g\|_\infty
\]
for every partition $P$ of $\ab.$ Consequently, the second inequality from \eqref{sum-2}
and estimate \eqref{est-2} are true also for the Young integral.
\end{proof}

\skipaline

Next convergence results are also true for all the three integrals under consideration.
For the Kurzweil-Stieltjes integral the proof is available e.g. in Chapter 6 of \cite{MST}.
The idea is pretty transparent and, as we will see below, applicable also to the Young and
Dushnik integrals: First, we notice that in both situations the sequences of integrals depending
on $n$ are Cauchy sequences in $\R$ and therefore they have a limit $I\in\R.$ Further, assumptions
on the convergence of functions involved, the estimates given in Lemma \ref{L2} and the existence
of the integrals $\intab f_n\,\dd g$ imply that the limit integrals exist and equals~$I.$

\skiphalfaline

\begin{The}\label{T3}
Let $f,g:[a,b]\to\R$ and $f_n:[a,b]\to\R,$ $n\in\N,$ be such that the integrals
$\intab f_n\,\dd g$ exist for all $n\in\N.$ Suppose at least one of the following
conditions is satisfied:
\begin{itemize}
\item[\rm(i)] $g\in\BV(\ab)$ and $f_n\rightrightarrows f$ on $\ab.$

\item[\rm(ii)] $g$ is bounded on $\ab$ and $\lim_{n\to\infty}\,\|f_n-f\|_{\BV}=0.$
\end{itemize}
Then the integral $\intab f\,\dd g$  exists as well, and
\[
   \lim_{n\to\infty}\intab f_n\,\dd g=\intab f\,\dd g.
\]
\end{The}
\begin{proof}
Both the integral and sum symbols now may refer to any of those three integrals we are
considering in this paper.

\skiphalfaline

We claim that in both cases (i) and (ii) the sequence $\{\intab f_n\,\dd g\}$ satisfies
the Cauchy condition. In case (i), we have by \eqref{est-1} from Lemma~\ref{L2}
\[
    \Big|\intab f_n\,\dd g-\intab f_m\,\dd g\Big|
   =\Big|\intab (f_n-f_m)\,\dd g\Big|\le\|f_n-f_m\|\,\varab\,g
\]
for all $m,n\in\N.$ Since $\varab\,g$ is finite and $\{f_n\}$ is uniformly convergent,
the right-hand side will be arbitrarily small if $m,n$ are sufficiently large.

\skiphalfaline

In case (ii), we use \eqref{est-2} from Lemma \ref{L2} to get
\begin{align*}
	 \Big|\intab f_n\,\dd g-\intab f_m\,\dd g\Big|
   &=\Big|\intab (f_n-f_m)\,\dd g\Big|\le 2\,\|g\|_{\infty}\,\|f_n-f_m\|_{\BV}
\\
   &\le 2\,\|g\|_{\infty}\,\|f_n-f\|_{\BV}+2\,\|g\|_{\infty}\,\|f-f_m\|_{\BV}
\end{align*}
for all $m,n\in\N.$ Since $\|g\|_\infty<\infty$ and $\|f_n-f\|_{\BV}\to 0,$ the right-hand
side will be arbitrarily small for $m,n$ sufficiently large.

\skiphalfaline

Thus, in both cases, there exists a~number $I\in\R$ such that
\begin{equation}\label{r6}
   \lim_{n\to\infty}\intab f_n\,\dd g=I.
\end{equation}

\skiphalfaline

To show that $\intab f\,\dd g=I,$ let $\eps>0$ be given. We claim there exists an $n_1\in\N$
such that
\begin{equation}\label{r7}
   |S(f-f_{n},\dd g,P)|<\eps\quad\mbox{for\ } n\ge n_1
                            \mbox{\ and every partition\ $P$ of $[a,b]$}.
\end{equation}

In case (i), this follows from \eqref{sum-1} in Lemma~\ref{L2}, which yields
\[
   |S(f-f_{n},\dd g,P)|\le\|f-f_n\|\,\varab\,g.
\]
In case (ii), we use \eqref{sum-2} of Lemma~\ref{L2}  to get
\[
   |S(f-f_{n},\dd g,P)|\le 2\|f-f_n\|_{\BV}\,\|g\|.
\]

These estimates show the validity of \eqref{r7}. By \eqref{r6}, there exists
an $n_0\ge n_1$ such that
\begin{equation*}
   \big|\intab f_{n_0}\,\dd g-I\big|<\eps,
\end{equation*}

Now, in the case of the Kurzweil-Stieltjes integral we can choose a gauge $\delta_\eps$
on $\ab$ such that
\begin{equation}\label{r8}
   |S(f_{n_0},\dd g,P)-\intab f_{n_0}\,\dd g|<\eps
\end{equation}
holds for each $\delta_\eps$-fine partition $P$ of $\ab.$

\skiphalfaline

Similarly, in the case of the Dushnik integral, we can choose a division $\balfa_0$ of $\ab$
such that \eqref{r8} holds for each partition $P=(\balfa,\ksi)$ of $\ab$ such that
$\balfa\supset\balfa_\eps$ and $\xi_j\in (\alpha_{j-1},\alpha_j)$ for $j\in\{1,\dots,\nu(P)\}.$

\skiphalfaline

Finally, in the case of the Young integral, we can choose a division $\balfa_\eps$ of $\ab$
such that
\[
   |S_Y(f_{n_0},\dd g,P)-\intab f_{n_0}\,\dd g|<\eps
\]
holds whenever $P=(\balfa,\ksi),$ $\balfa\supset\balfa_\eps$ and $\xi_j\in (\alpha_{j-1},\alpha_j)$
for all $j\in\{1,\dots,\nu(P)\}.$

\skiphalfaline

To summarize, in case of the Kurzweil-Stieltjes we have
\[
  |S(f,\dd g,P)-I|\le|S(f-f_{n_0},\dd g,P)|+|S(f_{n_0},\dd g,P)-I_{n_0}|+|I_{n_0}-I|
  <3\eps
\]
for each $\delta_\eps$-fine partition $P$ of $\ab,$ in case of the Dushnik integral
\[
  |S(f,\dd g,P)-I|<3\eps
\]
holds for each partition $P=(\balfa,\ksi)$ such that $\balfa\supset\balfa_\eps$ and
$\xi_j\in (\alpha_{j-1},\alpha_j)$  for all $j\in\{1,\dots,\nu(P)\},$ and in the case
of the Young integral
\[
  |S_Y(f,\dd g,P)-I|<3\,\eps
\]
holds for each partition $P=(\balfa,\ksi)$ such that $\balfa\supset\balfa_\eps$ and
$\xi_j\in (\alpha_{j-1},\alpha_j)$  for all $j\in\{1,\dots,\nu(P)\}.$ Thus, $\intab f\,\dd g=I$
holds in any of the considered cases. The proof is complete.
\end{proof}

\skipaline

Next convergence result is complementary to Theorem \ref{T3}. Its proof is based on
the same schedule as that of Theorem \ref{T3} and we leave it to the reader.

\skiphalfaline

\begin{The}\label{T4}
Let $f,g:[a,b]\to\R$ and $g_n:[a,b]\to\R,$ $n\in\N,$ be such that the integrals
$\intab f\,\dd g_n$ exist for all $n\in\N.$ Suppose at least one of the following
conditions is satisfied:
\begin{itemize}
\item[\rm(i)] $f\in\BV(\ab)$ and $g_n\rightrightarrows g$ on $\ab.$

\item[\rm(ii)] $f$ is bounded on $\ab$ and $\lim_{n\to\infty}\varab(g_n-g)=0.$
\end{itemize}
Then the integral $\intab f\,\dd g$  exists as well, and
\[
   \lim_{n\to\infty}\intab f\,\dd g_n=\intab f\,\dd g.
\]
\end{The}

\bigskip

Now we are able to prove Theorem \ref{main}.

\skiphalfaline

\noindent
{\bf Proof of Theorem \ref{main}}

\skiphalfaline

\noindent a) \ First, assume that $g\in\BV(\ab)$ and $f\in\G(\ab).$ Choose a sequence
$\{f_n\}$ of finite step functions such that $f_n\rightrightarrows f$ on $\ab.$
Then by Lemma \ref{L1} we have
\[
  {\rm(K)}\!\intab\,f_n\,\dd g={\rm(Y)}\!\intab\,f_n\,\dd g \quad\mbox{for all \ } n\in\N
\]
and further, by Theorem \ref{T3},
\[
  {\rm(K)}\!\intab\,f\,\dd g=\lim_{n\to\infty}{\rm(K)}\!\intab\,f_n\,\dd g
  =\lim_{n\to\infty}{\rm(Y)}\!\intab\,f_n\,\dd g={\rm(Y)}\!\intab\,f\,\dd g.
\]
For the Dushnik integral we have by Lemma \ref{L1} and Theorem \ref{T3}
\begin{align*}
  &{\rm(D)}\!\intab\,f\,\dd g=\lim_{n\to\infty}{\rm(D)}\!\intab\,f_n\,\dd g
  =\lim_{n\to\infty}\Big(f_n(b)\,g(b)-f_n(a)\,g(a)-{\rm(K)}\!\intab\,f_n\,\dd g\Big)
 \\
  &\quad=f(b)\,g(b)-f(a)\,g(a)-{\rm(K)}\!\intab\,f\,\dd g.
\end{align*}
Hence \eqref{eq-main} is true.

\skiphalfaline

\noindent b) \ It remains to consider the case that $f\in\BV(\ab)$ and $g\in\G(\ab).$
Choose a sequence $\{g_n\}$ of finite step functions such that $g_n\rightrightarrows g$
on $\ab.$ Then by Lemma \ref{L1} we have
\[
  {\rm(K)}\intab\,f\,\dd g_n={\rm(Y)}\intab\,f\,\dd g_n \quad\mbox{for all \ } n\in\N
\]
and, further, by Theorem \ref{T4},
\[
  {\rm(K)}\intab\,f\,\dd g=\lim_{n\to\infty}{\rm(K)}\intab\,f_n\,\dd g
  =\lim_{n\to\infty}{\rm(Y)}\intab\,f_n\,\dd g={\rm(Y)}\intab\,f\,\dd g.
\]
For the Dushnik integral we get using by Lemma \ref{L1} and Theorem \ref{T4}
\begin{align*}
  {\rm(D)}\intab\,f\,\dd g&=\lim_{n\to\infty}{\rm(D)}\intab\,f_n\,\dd g
 \\
  &=\lim_{n\to\infty}\Big(f_n(b)\,g(b)-f_n(a)\,g(a)-{\rm(K)}\intab\,f_n\,\dd g
 \\
  &=f(b)\,g(b)-f(a)\,g(a)-{\rm(K)}\intab\,f\,\dd g.
\end{align*}
This completes the proof.

\skipaline

\end{document}